\documentclass[twoside,12pt,leqno]{amsart}
\usepackage{hyperref}
\usepackage[all,cmtip]{xy}
\usepackage{amssymb}

\advance\topmargin by -\headheight\advance\topmargin by -\headsep
\numberwithin{equation}{section}
\newtheorem{theorem}{Theorem}[section]

\newtheorem{proposition}[theorem]{Proposition}

\theoremstyle{definition}
\newtheorem{definition}[theorem]{Definition}
\newtheorem{example}[theorem] {Example}

\newtheorem{remark}[theorem]{Remark}

\allowdisplaybreaks
\begin{document}
\renewcommand{\leq}{\leqslant}
\renewcommand{\le}{\leqslant}
\renewcommand{\geq}{\geqslant}
\renewcommand{\ge}{\geqslant}
\renewcommand{\setminus}{\smallsetminus}
\renewcommand{\emptyset}{\varnothing}
\renewcommand{\phi}{\varphi}
\newcommand{\suchthat}{\;\;|\;\;}
\newcommand{\abs}[1]{\lvert#1\rvert}
\newcommand{\norm}[1]{\lVert#1\rVert}
\newcommand{\la}{\langle}
\newcommand{\ra}{\rangle}
\newcommand{\into}{\hookrightarrow}
\newcommand{\lto}{\longrightarrow}
\newcommand{\lra}{\longrightarrow}
\newcommand{\xra}{\xrightarrow}

\newcommand{\dbar}{\bar{\partial}}

\newcommand{\R}{\mathbb{R}}
\newcommand{\V}{\mathbb{V}}
\newcommand{\Z}{\mathbb{Z}}
\newcommand{\C}{\mathbb{C}}
\newcommand{\E}{\mathbb{E}}
\newcommand{\HH}{\mathbb{H}}
\newcommand{\PP}{\mathbb{P}}

\newcommand{\cD}{\mathcal{D}}
\newcommand{\cE}{\mathcal{E}}
\newcommand{\cM}{\mathcal{M}}
\newcommand{\cO}{\mathcal{O}}
\newcommand{\cR}{\mathcal{R}}
\newcommand{\cS}{\mathcal{S}}

\newcommand{\PU}{\mathrm{PU}}
\newcommand{\PSU}{\mathrm{PU}}
\newcommand{\PGL}{\mathrm{PGL}}
\newcommand{\PSL}{\mathrm{PSL}}
\newcommand{\PSO}{\mathrm{PSO}}
\newcommand{\SU}{\mathrm{SU}}
\newcommand{\U}{\mathrm{U}}
\newcommand{\OO}{\mathrm{O}}
\newcommand{\GL}{\mathrm{GL}}
\newcommand{\SL}{\mathrm{SL}}
\newcommand{\SSS}{\mathrm{S}}
\newcommand{\SO}{\mathrm{SO}}
\newcommand{\Sp}{\mathrm{Sp}}
\newcommand{\Spin}{\mathrm{Spin}}
\newcommand{\Pin}{\mathrm{Pin}}
\newcommand{\PSpin}{\mathrm{PSpin}}

\newcommand{\Jac}{\operatorname{Jac}}
\newcommand{\ad}{\operatorname{ad}}
\newcommand{\Ad}{\operatorname{Ad}}
\newcommand{\tr}{\operatorname{tr}}
\newcommand{\rk}{\operatorname{rk}}
\newcommand{\im}{\operatorname{im}}
\newcommand{\coker}{\operatorname{coker}}
\newcommand{\Hom}{\operatorname{Hom}}
\newcommand{\End}{\operatorname{End}}
\newcommand{\Ext}{\operatorname{Ext}}
\newcommand{\Id}{\operatorname{Id}}
\newcommand{\Quot}{\operatorname{Quot}}
\newcommand{\codim}{\operatorname{codimension}}
\newcommand{\Index}{\operatorname{index}}
\newcommand{\aut}{\operatorname{aut}}
\newcommand{\Aut}{\operatorname{Aut}}
\newcommand{\SAut}{\operatorname{SAut}}
\newcommand{\Gr}{\operatorname{Gr}}
\newcommand{\Pic}{\operatorname{Pic}}

\hyphenation{Higgs}

\newcommand{\catq}{/\!\!/}

\newcommand{\liem}{\mathfrak{m}}
\newcommand{\liez}{\mathfrak{z}}
\newcommand{\liemp}{\mathfrak{m}_+}
\newcommand{\liemm}{\mathfrak{m}_-}
\newcommand{\liemc}{\mathfrak{m}^{\mathbb{C}}}
\newcommand{\lieh}{\mathfrak{h}}
\newcommand{\liehc}{\mathfrak{h}^{\mathbb{C}}}
\newcommand{\lieg}{\mathfrak{g}}
\newcommand{\liegc}{\mathfrak{g}^{\mathbb{C}}}
\newcommand{\vol}{\mathrm{vol}}

\newcommand{\lie}{\mathfrak}
\newcommand{\alie}{\mathfrak{a}}
\newcommand{\blie}{\mathfrak{b}}
\newcommand{\clie}{\mathfrak{c}}
\newcommand{\glie}{\mathfrak{g}}
\newcommand{\hlie}{\mathfrak{h}}
\newcommand{\klie}{\mathfrak{k}}
\newcommand{\llie}{\mathfrak{l}}
\newcommand{\mlie}{\mathfrak{m}}
\newcommand{\olie}{\mathfrak{o}}
\newcommand{\plie}{\mathfrak{p}}
\newcommand{\slie}{\mathfrak{s}}
\newcommand{\tlie}{\mathfrak{t}}
\newcommand{\ulie}{\mathfrak{u}}
\newcommand{\zlie}{\mathfrak{z}}
\newcommand{\sllie}{\mathfrak{sl}}
\newcommand{\solie}{\mathfrak{so}}
\newcommand{\splie}{\mathfrak{sp}}
\newcommand{\sulie}{\mathfrak{su}}
\newcommand{\gllie}{\mathfrak{gl}}


\markboth{Peter B. Gothen}{Hitchin pairs}
$ $
\bigskip

\bigskip

\centerline{{\Large Hitchin Pairs for non-compact real Lie groups}}

\bigskip
\bigskip
\centerline{{\large by Peter B.\ Gothen\footnote{The author was partially supported by CMUP (UID/MAT/00144/2013) and
the project PTDC/MAT-GEO/2823/2014 funded by FCT (Portugal) with
national and where applicable European structural funds through the
programme FEDER, under the partnership agreement PT2020.  The author
acknow\-ledges support from U.S. National Science Foundation grants
DMS 1107452, 1107263, 1107367 "RNMS: GEometric structures And
Representation varieties" (the GEAR Network)}
}}

\vspace*{.7cm}

\noindent
  \begin{changemargin}{1cm}{1cm}
    \small
    \begin{center}
      \textbf{Abstract}
    \end{center}
    Hitchin pairs on Riemann surfaces are generalizations of Higgs bundles, allowing the Higgs field to be twisted by an arbitrary line bundle. We consider this generalization in the context of $G$-Higgs bundles for a real reductive Lie group $G$.  We outline the basic theory and review some selected results, including recent results by Nozad and the author \cite{gothen-nozad:2016} on Hitchin pairs for the unitary group of indefinite signature $\U(p,q)$.
  \end{changemargin}

\vspace*{.7cm}

\pagestyle{myheadings}
\section{Introduction}
\label{sec:introduction}

Let $X$ be a closed Riemann surface with holomorphic cotangent bundle
$K= \Omega^1_X$.  A rank $n$ \emph{Higgs bundle} on $X$ is a pair
$(E,\phi)$, where $E\to X$ is a rank $n$ holomorphic vector bundle and
$\phi\colon E \to E\otimes K$ is an endomorphism valued holomorphic
$1$-form on $X$. Higgs bundles are fundamental objects in the
non-abelian Hodge theorem
\cite{corlette:1988,donaldson:1987,hitchin:1987a,simpson:1988}. In the
simplest (abelian) case of $n=1$ this can be expressed as the
isomorphism
\begin{displaymath}
   \Hom(\pi_1X,\C^*) \simeq  T^*\Jac(X),
\end{displaymath}
whose infinitesimal version 
gives the Hodge decomposition 
\begin{math}
  H^1(X,\C) \simeq H^{1,0}(X) \oplus H^{0,1}(X).
\end{math}
Thus, for $n=1$, a flat line bundle on $X$ corresponds to a pair
$(E,\phi)$ consisting of a holomorphic line bundle $E\to X$ and a
holomorphic $1$-form $\phi$ on $X$. For general $n$, non-abelian Hodge
theory produces an isomorphism
\begin{displaymath}
  \Hom(\pi_1X,\GL(n,\C))\catq \GL(n,\C) \simeq 
  \mathcal{M}(\GL(n,\C)). 
\end{displaymath}
Here the space on the right hand side is the moduli space of
isomorphism classes of Higgs
bundles (of degree 0) and the space on the left hand side is
the space of representations of $\pi_1X$ modulo the action of
$\GL(n,\C)$ by overall conjugation.  
Note that, viewed in this way, the non-abelian Hodge theorem
generalizes the Narasimhan--Seshadri theorem \cite{narasimhan-seshadri:1965}
to non-compact groups.

For many purposes, rather than considering $\phi$ as a $1$-form, one
might as well consider pairs $(E,\phi)$, where $\phi\colon E\to
E\otimes L$ is twisted by an arbitrary line bundle $L \to X$. Such a
pair is known as a \emph{Hitchin pair} or a \emph{twisted Higgs
  bundle}. This point of view was probably first explored
systematically by Nitsure \cite{nitsure:1991}.  The non-abelian Hodge
theorem generalizes to this context and involves, on the one side,
meromorphic Higgs bundles and on the other side meromorphic
connections. This generalization has been carried out by Simpson
\cite{simpson:1990}, for Higgs fields with simple poles, and
Biquard--Boalch \cite{biquard-boalch:2004}, for more general polar parts
(see Boalch \cite{boalch:2012} for a survey).

Another generalization of the non-abelian Hodge theorem has to do with
representations of $\pi_1X$ in groups $G$ other than the general
linear group. This already goes back to Hitchin's seminal papers
\cite{hitchin:1987a,hitchin:1987b,hitchin:1992} and indeed was also
treated by Simpson \cite{simpson:1992}. Here we shall focus on the
theory for real $G$, which has quite a different flavour from the
theory for complex $G$.  A systematic approach to non-abelian Hodge
theory for real reductive groups $G$ and applications to the study of
character varieties has been explored in a number of papers; see, for
example,
\cite{gothen:2001,bradlow-garcia-prada-gothen:hss-higgs,garcia-gothen-mundet:2009b,garcia-gothen-mundet:2013}.
The focus of the present paper are the objects which are obtained by
allowing for an arbitrary twisting line bundle $L$ in $G$-Higgs
bundles rather than just the canonical bundle $K$. These objects are
known as \emph{$G$-Hitchin pairs}.

There are many other important aspects of Higgs bundle theory and
without any pretense of completenes, we mention here a few.  One of
the important features of the Higgs bundle moduli space for complex
$G$ is that it is an algebraically completely integrable Hamiltonian
system (see Hitchin\cite{hitchin:1987b}), known as the Hitchin
system. This is closely related to the fact that this moduli space is
a holomorphic symplectic manifold admitting a hyper-K\"ahler metric.
This aspect of the theory can be generalized to Hitchin pairs using
Poisson geometry, as pioneered by Bottacin \cite{bottacin:1995} and
Markman \cite{markman:1994}; see Biquard--Boalch
\cite{biquard-boalch:2004} for the existence of hyper-K\"ahler metrics
on the symplectic leaves.  Closely related is the theory of
parabolic Higgs bundles (see, for example, Konno \cite{konno:1993} and
Yokogawa \cite{yokogawa:1993}). Parabolic $G$-Higgs bundles for real
$G$ have been considered by, among others,
Logares \cite{Logares:2006}, Garc{\'{\i}}a-Prada--Logares--Mu{\~n}oz
\cite{garcia-logares-munoz:2009} and
Biquard--Garc{\'{\i}}a-Prada--Mundet
\cite{biquard-garcia-mundet:2015}. Higgs bundles also play an
important role in mirror symmetry (see, for example, Hausel--Thaddeus
\cite{hausel-thaddeus:2003}) and in the geometric Langlands
correspondence (see, for example, Kapustin--Witten
\cite{kapustin-witten:2007}).  Also a number of results on $G$-Higgs
bundles for real groups can be obtained via the study of the Hitchin
fibration; for this we refer the reader to
Baraglia--Schaposnik \cite{baraglia-schaposnik:2015}, Garc{\'\i}a-Prada--Pe\'on-Nieto--Ramanan
\cite{garcia-peon-ramanan:2015}, Hitchin--Schaposnik
\cite{hitchin-schaposnik:2014}, Pe\'on-Nieto \cite{peon:2015} and
Schaposnik \cite{schaposnik:2015}, as well as further references found
therein.


In this paper we describe the basics of the theory of $G$-Hitchin
pairs and give a few examples
(Section~\ref{sec:hitchin-pairs-real}). We explain the
Hitchin--Kobayashi correspondence which relates the (parameter
dependent) stability condition for $G$-Hitchin pairs to solutions to
Hitchin's gauge theoretic equations (Section~\ref{sec:hk}). 
We then describe recent work of Nozad
and the author \cite{gothen-nozad:2016} on $\U(p,q)$-Hitchin pairs
(introduced in Section~\ref{sec:upq-hp}), the Milnor--Wood inequality
for such pairs (Section~\ref{sec:mw}) and how wall-crossing arguments
can be used to study their moduli (Section~\ref{sec:wall-crossing}).

\subsubsection*{Acknowledgments}
I would like to thank P.~Boalch for useful comments on the first
version of the paper.

This paper is based on my talk at the GEOQUANT 2015 conference held at
the Instituto de Ciencias Matem\'aticas in Madrid. It is a great
pleasure for me to thank the organizers for the invitation to speak
and for organizing such an enjoyable and stimulating meeting.

\section{Hitchin pairs for real groups}
\label{sec:hitchin-pairs-real}

Let $G$ be a connected real reductive Lie group. Following Knapp
\cite{knapp:1996}, we shall take this to mean that the following data
has been fixed:
\begin{itemize}
\item a maximal compact subgroup $H \subset G$;
\item a Cartan decomposition $\lieg = \lieh + \liem$;
\item a non-degenerate $\Ad(G)$-invariant quadratic form, negative
  definite on $\lieh$ and positive definite on $\liem$, which
  restricts to the Killing form on the semisimple part
  $\lieg_{\mathrm{ss}} = \lieg / [\lieg,\lieg]$ of $\lieg$.
\end{itemize}
Note that the above data complexify (with the possible exception of
$G$) and that there is an isotropy representation
\begin{displaymath}
  \iota\colon H^\C\to \Aut(\liemc)
\end{displaymath}
coming from restricting and complexifying the adjoint representation
of $G$.

Let $X$ be Riemann surface and let $K = \Omega^1_X$ be its holomorphic
cotangent bundle. Fix a line bundle $L\to X$. For a principal
$H^\C$-bundle $E\to X$ and a representation $\rho\colon H^\C\to\GL(V)$
of $H^\C$, we denote the associated vector bundle by $E(V) =
E\times_\rho V$.

\begin{definition}
  A \emph{$G$-Hitchin pair} (twisted by $L$) on $X$ is a pair
  $(E,\phi)$, where $E\to X$ is a holomorphic principal $H^\C$-bundle
  and $\phi\in H^0(X,L\otimes E(\liemc))$ is a holomorphic $1$-form
  with values in the vector bundle defined by the isotropy
  representation of $H^\C$. If $L=K$, the pair $(E,\phi)$ is called a
  \emph{$G$-Higgs bundle}.
\end{definition}

\begin{example}
  If $G$ is compact, a $G$-Hitchin pair is nothing but a holomorphic
  principal $G^\C$-bundle.
\end{example}

\begin{example}
  If $G=\GL(n,\C)$, a $G$-Hitchin pair is a pair $(E,\phi)$, where
  $E\to X$ is a rank $n$ holomorphic vector bundle and $\phi\in
  H^0(X,L\otimes\End(E))$ is an $L$-twisted endomorphism of $E$. A
  $\SL(n,\C)$-Higgs bundle is given by the same data, with the
  additional requirements that $\det(E)=\mathcal{O}_X$ and $\phi\in
  H^0(X,L\otimes\End_0(E))$, where $\End_0(E)\subset\End(E)$ is the
  subbundle of $\phi$ with $\phi=0$.
\end{example}

\begin{example}
  Let $G=\SL(n,\R)$. A maximal compact subgroup is $\SO(n)$ defined
  by the standard inner product $\la x,y \ra = \sum x_iy_i$ and the
  isotropy representation is the subspace of $A\in\mathfrak{sl}(n,R)$
  which are symmetric with respect to the inner product:
  \begin{displaymath}
    \la Ax,y \ra = \la x,Ay \ra.
  \end{displaymath}
  Hence a $\SL(n,\R)$-Hitchin pair can be viewed as a pair
  $((U,Q),\phi)$, where $(U,Q)$ is a holomorphic orthogonal bundle,
  i.e., $U\to X$ is a rank $n$ vector bundle with a non-degenerate
  holomorphic quadratic form $Q$, and $\phi\in H^0(X,L\otimes
  S^2_QU)$. Here $S^2_QU\subset \End(U)$ denotes the subbundle of
  endomorphism of $U$, which are symmetric with respect to $Q$.
\end{example}

\begin{example}
  \label{exa:upq}
  Let $G=\U(p,q)$, the group of linear transformations of $\C^{p+q}$
  which preserves an indefinite hermitian form of signature $(p,q)$ on
  $\C^{p+q} = \C^p\times\C^q$. Taking the obvious $\U(p)\times\U(q)$
  as the maximal compact subgroup, we have
  $H^\C=\GL(p,\C)\times\GL(q,\C)$ and the isotropy representation is
  \begin{displaymath}
    \GL(p,\C)\times\GL(q,\C) \to
    \Hom(\C^q,\C^p)\oplus\Hom(\C^p,\C^q)    
  \end{displaymath}
  acting by restricting the adjoint representation of $\GL(p+q,\C)$.
  Hence a $\U(p,q)$-Hitchin pair can be identified with a quadruple
  $(V,W,\beta,\gamma)$, where
  \begin{displaymath}
    \beta\in H^0(L\otimes\Hom(W,V))
    \quad\text{and}\quad
    \gamma\in H^0(L\otimes\Hom(V,W)).
  \end{displaymath}
  The $\GL(p+q,\C)$-Hitchin pair associated via the inclusion
  $\U(p,q)\subset\GL(p+q,\C)$ is $(E,\phi)$, where $E=V\oplus W$ and
  $\phi=\left(
    \begin{smallmatrix}
      0 & \beta \\
      \gamma & 0
    \end{smallmatrix}
  \right)$.
  Of course a $\SU(p,q)$-Hitchin pair is given by the same data, with
  the additional requirement that
  $\det(V)\otimes\det(W)=\mathcal{O}_X$. 
\end{example}

\begin{example}
  Let $G=\Sp(2n,\R)$, the real symplectic group in dimension $2n$,
  defined as the subgroup of $\SL(2n,\R)$ of transformations of
  $\R^{2n}$ preserving the standard symplectic form, which can be
  written in coordinates $(x_1,y_1,\dots,x_n,y_n)\in\R^ {2n}$ as
  \begin{displaymath}
    \omega = dx_1\wedge dy_1 + \dots + dx_n\wedge dy_n.
  \end{displaymath}
  Then a $\Sp(2n,\R)$-Hitchin pair can identified with a triple
  $(V,\beta,\gamma)$, where $V\to X$ is a rank $n$ vector bundle and
    \begin{displaymath}
    \beta\in H^0(L\otimes S^2V)
    \quad\text{and}\quad
    \gamma\in H^0(L\otimes S^2V^*).
  \end{displaymath}
  Note how the inclusions $\Sp(2n,\R)\subset\SL(2n,\R)$ and
  $\Sp(2n,\R)\subset\SU(n,n)$ are reflected in the
  associated vector bundle data. In the former case, the rank $2n$
  orthogonal bundle $(U,Q)$ is given by $U=V\oplus V^*$ with the
  quadratic form
  $Q=\left(
    \begin{smallmatrix}
      0 & 1 \\
      1 & 0
    \end{smallmatrix}
  \right)$.
\end{example}

\section{The Hitchin--Kobayashi correspon\-dence}
\label{sec:hk}

We now move on to the central notion of stability for $G$-Hitchin
pairs.  The stability condition depends on a parameter $c\in
i\liez$, where $\liez$ denotes the centre of $\lieh$. 

From the point of view of construction of moduli spaces, stability
allows for a GIT construction of the moduli space
$\mathcal{M}_d^c(X,G)$ of $c$-semistable $G$-Higgs bundles for a fixed
topological invariant $d\in \pi_1(H)$; this construction has been
carried out by Schmitt (see \cite{schmitt:2008}).

On the other hand, there is a Hitchin--Kobayashi correspondence for
$G$-Higgs bundles, which gives necessary and sufficient conditions in
terms of stability for the existence of solutions to the so-called
Hitchin's equations. To state these equations, we need some notation.
By a hermitian metric on the $H^\C$-bundle $E$ we mean a reduction of
structure group to $H\subset H^\C$, i.e., a smooth section $h\colon X
\to E(H^\C/H)$. We denote the corresponding principal $H$-bundle by
$E_h$. Note that $h$ defines a compact real structure, denoted by
$\sigma_h$, on the bundle of Lie algebras $E(\liegc)$, compatible with
the decomposition $E(\liegc)=E(\liehc)\oplus E(\liemc)$. If we combine
$\sigma_h$ with the conjugation on complex $1$-forms on $X$, we obtain a
complex antilinear involution
\begin{math}
  A^1(E(\liegc))\to A^1(E(\liegc)).
\end{math}
This restricts to an antilinear map which, by a slight abuse of
notation, we denote by the same symbol:
\begin{displaymath}
  \sigma_h\colon A^{1,0}(E(\liemc)) \to A^{0,1}(E(\liemc)).
\end{displaymath}
Fix a hermitian metric $h_L$ on $L$ and let $\omega_X$ denote the
K\"ahler form of a metric on $X$ compatible with its complex
structure, normalized so that $\int_X\omega_X=2\pi$.
Then,
for $c\in i\liez$, Hitchin's equation for a metric $h$ on $E$ is
the following
\begin{equation}
  \label{eq:hitchin}
  F(A_h) + [\phi,\sigma_h(\phi)]\omega_X = -ic\omega_X. 
\end{equation}
Here $A_h$ denotes the Chern connection on $E_h$ (i.e., the unique
$H$-connection compatible with the holomorphic structure on $E$) and
$F(A_h)$ its curvature. Moreover, the bracket $[\phi,\sigma_h(\phi)]$
is defined by combining the Lie bracket on $\liegc = \liehc + \liemc$
with the contraction $L\otimes\overline{L}\to\mathcal{O}_X$ given by
the metric $h_L$. Note also that in the case when $L = K$, the second
term on the left hand side can be written simply as
\begin{math}
  [\phi,\sigma_h(\phi)]
\end{math}
where the bracket on the Lie algebra is now combined with the wedge
product on forms.

In order to state the Hitchin--Kobayashi correspondence for
$G$-Hitchin pairs, giving necessary and sufficient conditions for the
existence of solutions to the Hitchin equation, one needs an
appropriate stability condition. The general condition needed can be
found in \cite{garcia-gothen-mundet:2009b} (based, in turn, on
Bradlow--Garc{\'\i}a-Prada-Mundet
\cite{bradlow-garcia-prada-mundet:2003} and Mundet
\cite{mundet:2000}). It is fairly involved to state in general, so we
shall refer the reader to loc.\ cit.\ for the full statement and here
just give a couple of examples which cover our present needs. Note
that, just as the Hitchin equation, the stability condition will
depend on a parameter $c\in i\liez$. 

\begin{example}
  \label{exa:glnc-stability} (Cf.\ Hitchin\cite{hitchin:1987a}, Simpson \cite{simpson:1988,simpson:1992}.)
  Consider $\GL(n,\C)$-Hitchin pairs $(E,\phi)$, where $E\to X$ is a
  rank $n$ vector bundle and $\phi\in H^0(X,L\otimes\End(E))$. Recall
  that the \emph{slope} of a vector bundle $E$ on $X$ is the ratio
  between its degree and its rank: $\mu(E) = \deg(E) / \rk(E)$. A
  $\GL(n,\C)$-Hitchin pair $(E,\phi)$ is \emph{semistable} if
  \begin{equation}
    \label{eq:slope}
    \mu(F) \leq \mu(E)
  \end{equation}
  for all non-zero subbundles $F\subset E$ which are preserved by
  $\phi$, i.e., such that $\phi(F)\subset F\otimes L$. Moreover,
  $(E,\phi)$ is \emph{stable} if additionally strict inequality holds
  in (\ref{eq:slope}) whenever $F\neq E$. Finally, $(E,\phi)$ is
  polystable if it is the direct sum of stable Higgs bundles, all of
  the same slope. In this case $i\liez \simeq\R$ and the stability
  parameter is fixed to be the real constant $c=\mu(E)$. Note that
  this constraint is of a topological nature and can be obtained from
  Chern--Weil theory by integrating the trace of the Hitchin equation,
  which in this case is:
  \begin{displaymath}
    F(A_h) + [\phi,\phi^{\ast_h}]\omega_X = -ic\Id\omega_X.
  \end{displaymath}

\end{example}

\begin{example}(Cf.\ \cite{bradlow-garcia-prada-gothen:2003}.)
  \label{exa:upq-stability}
  Consider $\U(p,q)$-Hitchin pairs $(V,W,\beta,\gamma)$. In this case,
  $i\liez\simeq \R\times\R$ and the Hitchin equation becomes
  \begin{equation}
  \label{eq:upq-hitchin}
  \begin{aligned}
    F(A_h(V))+(\beta\beta^{\ast_h}
      -\gamma^{\ast_h}\gamma)\omega_X&=-ic_1\Id_{V}\omega_X,\\
    F(A_h(W))+(\gamma\gamma^{\ast_h}
      -\beta^{\ast_h}\beta)\omega_X&=-ic_2\Id_{W}\omega_X.
  \end{aligned}
  \end{equation}
  Here $A_h(V)$ and $A_h(W)$ denote the Chern connections on $V$ and
  $W$, respectively, and the parameter $(c_1,c_2)\in\R\times\R$ is
  constrained by Chern--Weil theory by
  \begin{displaymath}
    \frac{p}{p+q}c_1 +  \frac{q}{p+q}c_2 = \mu(V\oplus W).
  \end{displaymath}
  The stability condition is most conveniently described by introducing the
  \emph{$\alpha$-slope} of $(V,W,\beta,\gamma)$ by
  \begin{displaymath}
    \mu_\alpha(V,W,\beta,\gamma) = \mu(V\oplus W) + \alpha\frac{p}{p+q}
  \end{displaymath}
  for a real parameter $\alpha$, related to $(c_1,c_2)$ by
  $\alpha=c_2-c_1$. The $\alpha$-stability conditions are completely
  analogous to the ones of Example~\ref{exa:glnc-stability}, but
  applied to $\U(p',q')$-subbundles, defined in the obvious way by
  $V'\subset V$ and $W'\subset W$ such that $\beta(W')\subset
  V'\otimes L$ and $\gamma(V')\subset W'\otimes L$.
\end{example}

The Hitchin--Kobayashi correspondence for $G$-Hitchin pairs
\cite{hitchin:1987a,simpson:1992,bradlow-garcia-prada-mundet:2003,garcia-gothen-mundet:2009b}
can now be stated as follows.

\begin{theorem}
  \label{thm:HK}
  Let $(E,\phi)$ be a $G$-Hitchin pair. There exists a hermitian
  metric $h$ in $E$ solving Hitchin's equation (\ref{eq:hitchin}) if
  and only if $(E,\phi)$ is $c$-polystable. Moreover, the
  solution $h$ is unique up to $H$-gauge transformations of $E_h$.
\end{theorem}

Next we explain how to give an interpretation in terms of moduli
spaces. Fix a $C^\infty$ principal $H$-bundle $\mathcal{E}$ of
topological class $d\in\pi_1H$ and consider the configuration space of
$G$-Higgs pairs on $\mathcal{E}$:
\begin{displaymath}
  \mathcal{C}(\mathcal{E}) = \{(\dbar_A,\phi) \suchthat \dbar_A\phi=0\}.
\end{displaymath}
Here $\dbar_A$ is a $\dbar$-operator on $\mathcal{E}$ defining a
structure of holomorphic principal $H^\C$-bundle $E_A\to X$ and the
$C^\infty$-Higgs field $\phi\in A^{1,0}(\mathcal{E}(\liemc))$.  Let
$\mathcal{C}^{c-\mathrm{ps}}(\mathcal{E})\subset\mathcal{C}(\mathcal{E})$
be the subset of $c$-polystable $G$-Higgs pairs.  The \emph{complex gauge
group} $\mathcal{G}^\C$ is the group of $C^\infty$ automorphisms of the
principal $H^\C$-bundle $\mathcal{E}_{\C}$ obtained by extending the
structure group to the complexification $H^\C$ of $H$. It acts on
$\mathcal{C}^{c-\mathrm{ps}}(\mathcal{E})$ and we can identify, as
sets\footnote{Indeed a construction of the moduli space using complex analytic
methods in the style of Kuranishi should be possible, though we are
not aware of the existence of such a construction in the literature.},
\begin{displaymath}
  \mathcal{M}_d^c(X,G)
  = \mathcal{C}^{c-\mathrm{ps}}(\mathcal{E}) / \mathcal{G}^\C.
\end{displaymath}

Now consider Hitchin's equation (\ref{eq:hitchin}) as an equation for
a pair $(A,\phi)$ of a (metric) connection $A$ on $\mathcal{E}$ and a
Higgs field $\phi\in A^{1,0}(\mathcal{E}(\liemc))$.  The complex gauge
group $\mathcal{G}^\C$ acts transitively on the space of metrics on
$\mathcal{E}$ with stabilizer the \emph{unitary gauge group}
$\mathcal{G}$, by which we understand the $C^\infty$ automorphism
group of the $H$-bundle $\mathcal{E}$. Thus the Hitchin--Kobayashi
correspondence of Theorem~\ref{thm:HK} says that there is a complex
gauge transformation taking $(A,\phi)$ to a solution to Hitchin's
equation if and only if $(E_A,\phi)$ is $c$-polystable, and this
solution is unique up to unitary gauge transformation. In other words,
we have a bijection
\begin{equation}
  \label{eq:HK-moduli}
  \mathcal{M}^c_d(X,G) \simeq 
  \{(A,\phi)\suchthat\text{$(A,\phi)$ satisfies (\ref{eq:hitchin})}\}
  / \mathcal{G}.
\end{equation}

When $G$ is compact, there is no Higgs field and the Hitchin equation
simply says that the Chern connection is (projectively) flat. Hence
\eqref{eq:HK-moduli} identifies the moduli space of semistable $G^\C$-bundles
with the moduli space of (projectively) flat $G$-connections. This
latter space can in turn be identified with the \emph{character
  variety} of representations of (a central extension of) the
fundamental group of $X$ in $G$.

For non-compact $G$, assume that $L=K$ and that the parameter $c\in
iZ(\lieg)$. Then the Hitchin equation can be interpreted as a
(projective) flatness condition for the $G$-connection $B$ defined
by
\begin{equation}
  \label{eq:flat-G}
  B = A_h + \phi - \sigma_h(\phi).
\end{equation}
It is a fundamental theorem of Donaldson \cite{donaldson:1987} and
(more generally) Corlette \cite{corlette:1988} that for any flat
reductive\footnote{When $G$ is linear this simply means that the
  holonomy representation is completely reducible.} connection $B$ on
a principal $G$-bundle $\mathcal{E}_G$, there exists a so-called
{harmonic metric} on $\mathcal{E}_G$. A consequence of harmonicity is
that when the metric is used to decompose $B$ as in \eqref{eq:flat-G},
then $(A,\phi)$ satisfies the Hitchin equation. Combining this with
the Hitchin--Kobayashi correspondence gives the non-abelian Hodge
theorem\footnote{See the references cited in the Introduction for the
  generalization to the meromorphic situation.}: an identification
between the moduli space of $G$-Higgs bundles and the character
variety for representations of (a central extension of) $\pi_1 X$ in
$G$.

\begin{example}
  If we want to apply the non-abelian Hodge theorem to $\U(p,q)$-Higgs
  bundles, we need to fix the parameter in Hitchin's equation to be in
  the centre of $\U(p,q)$, i.e., in the notation of
  Example~\ref{exa:upq-stability}, we must take $c_1=c_2=c=\mu(V\oplus
  W)$. Of course this corresponds to the value for $\GL(n,\C)$-Higgs
  bundles under the inclusion $\U(p,q)\subset \GL(p+q,\C)$ (cf.\
  Examples \ref{exa:upq} and \ref{exa:glnc-stability}).
\end{example}

\section{Hitchin pairs for $\U(p,q)$ and quiver bundles}
\label{sec:upq-hp}

We saw in Example~\ref{exa:upq-stability}, that there is a degree of
freedom in the choice of stability parameter for $\U(p,q)$-Hitchin
pairs.  There is another way of viewing this parameter dependence for
the stability condition, which is to notice that a $\U(p,q)$-Hitchin
pair can be viewed as a \emph{quiver bundle} (see, e.g., King
\cite{king:1994}, \'Alvarez-C\'onsul--Garc{\'\i}a-Prada
\cite{alvarez-consul-garcia-prada:2003a,alvarez-consul-garcia-prada:2003b},
and also \cite{gothen-king:2005}). To explain this, recall that a
quiver $Q$ is an oriented graph (which we shall assume to be finite),
given by a set of vertices $Q_0$, a set of arrows $Q_1$ and head and
tail maps
\begin{displaymath}
  h,t\colon Q_1\to Q_0.
\end{displaymath}
For each $a\in Q_1$, let $M_a\to X$ be a holomorphic vector bundle on
$X$ and let $M=\{M_a\}$ be the collection of these \emph{twisting
  bundles}.

\begin{definition}
  A \emph{$Q$-bundle twisted by $M$} on $X$ is a collection of
  holomorphic vector bundles $E_i\to X$ indexed by the vertices $i\in
  Q_0$ of $Q$ and a collection of holomorphic maps $\phi_a\colon
  M_a\otimes E_{ta}\to E_{ha}$ indexed by the arrows $a\in Q_1$ of
  $Q$.
\end{definition}

\begin{remark}
  It is easy to see that $Q$-bundles on $X$ form a category which can
  be made into an abelian category by considering coherent
  $Q$-sheaves, in a way analogous to what happens for vector bundles. 
\end{remark}

It should now be clear that $L$-twisted $\U(p,q)$-Hitchin pairs  can be
viewed as $Q$-bundles for the quiver
\begin{equation}
  \label{eq:upq-quiver}
\xymatrix{
\bullet\ar@{<-}@/_1.2pc/[r]&\bullet\ar@{<-}@/_1.2pc/[l].
}
\medskip
\end{equation}
where both arrows are twisted by $L^*$.

There is a natural stability condition for quiver bundles which, just
as for Hitchin pairs, gives necessary and sufficient conditions for
the existence of solutions to certain natural gauge theoretic
equations (cf.\ King \cite{king:1994} and
\'Alvarez-C\'onsul--Garc{\'\i}a-Prada
\cite{alvarez-consul-garcia-prada:2003a,alvarez-consul-garcia-prada:2003b}).
This condition depends on a parameter vector
\begin{displaymath}
  \alpha = (\alpha_i)_{i\in Q_0} \in \R^{Q_0}
\end{displaymath}
and it is defined using the \emph{$\alpha$-slope} of a $Q$-bundle $E$:
\begin{displaymath}
  \mu_\alpha(E) = \frac{\sum_i(\deg(E_i) + \alpha_i\rk(E_i))}{\sum_i\rk(E_i)}. 
\end{displaymath}
Thus $E$ is $\alpha$-stable if for any proper non-zero sub-$Q$-bundle
$E'$ of $E$, we have
\begin{displaymath}
  \mu_\alpha(E') < \mu_\alpha(E), 
\end{displaymath}
and $\alpha$-semi- and polystability are defined just as for vector
bundles. 

Note that the stability condition is unchanged under an overall
translation of the stability parameter
\begin{displaymath}
  (\alpha_i) \mapsto (\alpha_i+a)
\end{displaymath}
for any constant $a\in\R$. Thus we may as well take $\alpha_0=0$ and
we see that the number of effective stability parameters is
$\abs{Q_0}-1$. In the case of $Q$-bundles for the quiver
\eqref{eq:upq-quiver}, i.e., $\U(p,q)$-Hitchin pairs, we then have one
real parameter $\alpha=\alpha_1$ and the general $Q$-bundle stability
condition reproduces the stability for $\U(p,q)$-Hitchin pairs of
Example \ref{exa:upq-stability}.

\section{The Milnor--Wood inequality for $\U(p,q)$-Hitchin pairs}
\label{sec:mw}

The Milnor--Wood inequality has its origins
\cite{milnor:1957,wood:1971} in the theory of flat
bundles. From this point of view there is a long sequence of
generalizations and important contributions (see, for example,
Dupont \cite{dupont:1978}, Toledo \cite{toledo:1979}, Domic--Toledo
\cite{domic-toledo:1987}, Turaev \cite{turaev:1984}, Clerc--{\O}rsted
\cite{clerc-orsted:2003}, Burger--Iozzi-Wienhard
\cite{burger-iozzi-wienhard:2003,burger-iozzi-wienhard:2010}). Here we
shall, however, focus on its Higgs bundle incarnation, again first
considered by Hitchin \cite{hitchin:1987a}. From this point of view it
is a bound on the topological class of a 
$\U(p,q)$-Hitchin pair. In order to state it we need the following
definition. 

\begin{definition}
  \label{def:toledo-invariant}
  Let $E=(V,W,\beta,\gamma)$ be a $\U(p,q)$-Hitchin pair. The
  \emph{Toledo invariant} of $E$ is
  \begin{displaymath}
    \tau(E) = \frac{2pq}{p+q}\big(\mu(V)-\mu(W)\big).
  \end{displaymath}
\end{definition}
Note that, if we set $a=\deg(V)$ and $b=\deg(W)$, then we can write
$\tau(E) = 2{(qa-pb)}/{(p+q)}$.

The Milnor--Wood inequality for $\U(p,q)$-Hitchin pairs can now be
stated as follows:

\begin{proposition}[Gothen--Nozad {\cite[Proposition~3.3]{gothen-nozad:2016}}]
  \label{prop:MW-inequality}
  Let $E=(V,W,\beta,\gamma)$ be an $\alpha$-semistable
  $\U(p,q)$-Hitchin pair with twisting line bundle $L$. Then
  \begin{displaymath}
    -\rk(\beta)\deg(L)+\alpha\big(\rk(\beta)-\frac{2pq}{p+q}\big)
    \leq \tau(E)\leq
    \rk(\gamma)\deg(L)+\alpha\big(\rk(\gamma)-\frac{2pq}{p+q}\big).
  \end{displaymath}
\end{proposition}

The proof is analogous to the one for $\U(p,q)$-Higgs bundles in
\cite{bradlow-garcia-prada-gothen:2003}. It applies the
$\alpha$-semistability condition for $\U(p,q)$-Hitchin pairs to
certain subobjects defined in a natural way using $\beta$ and
$\gamma$. We refer the reader to \cite{gothen-nozad:2016} for details.

\begin{remark}
  The Toledo invariant has been defined for $G$-Higgs bundles for any
  non-compact simple reductive group $G$ of hermitian type by
  Biquard--Garc{\'\i}a-Prada--Rubio
  \cite{biquard-garcia-rubio:2015}. These authors also prove a very
  general Milnor--Wood inequality for such $G$-Higgs bundles. In the
  case when $L=K$ their theorem specializes to our
  Proposition~\ref{prop:MW-inequality}.
\end{remark}

The inequality of Proposition~\ref{prop:MW-inequality} has several
interesting consequences, for example we get the following bounds on
the Toledo invariant (cf.\
\cite[Proposition~3.4]{gothen-nozad:2016}).

\begin{proposition}
  Let $E=(V,W,\beta,\gamma)$ be an $\alpha$-semistable
  $\U(p,q)$-Hitchin pair with twisting line bundle $L$ with
  $\deg(L)\geq 0$. Then the following hold:
  \begin{itemize}
  \item[$(i)$] If $\alpha \leq -\deg(L)$ then
    \begin{displaymath}
      \min\{p,q\}\bigl(-\alpha\frac{\abs{p-q}}{p+q}-\deg(L)\bigr)
      \leq \tau(E) \leq
      -\alpha\frac{2pq}{p+q}.
    \end{displaymath}
  \item[$(ii)$] If $-\deg(L)\leq\alpha\leq\deg(L)$ then
    \begin{displaymath}
      \min\{p,q\}\bigl(-\alpha\frac{\abs{p-q}}{p+q}-\deg(L)\bigr)
      \leq \tau(E) \leq
      \min\{p,q\}\bigl(\deg(L)-\alpha\frac{\abs{p-q}}{p+q}\bigr).
    \end{displaymath}
  \item[$(iii)$] If $\deg(L)\leq\alpha$ then
    \begin{displaymath}
      -\alpha\frac{2pq}{p+q}
      \leq \tau(E) \leq
      \min\{p,q\}\bigl(\deg(L)-\alpha\frac{\abs{p-q}}{p+q}\bigr).
    \end{displaymath}    
  \end{itemize}
\end{proposition}

Note, in particular, that for $\alpha = 0$ (the value relevant for the
non-abelian Hodge theorem) we have by $(ii)$ of the proposition that
\begin{equation}
  \label{eq:MW-L}
  \abs{\tau(E)} \leq\min\{p,q\}\deg(L).
\end{equation}
In the case of $\U(p,q)$-Higgs bundles (i.e., $\alpha=0$ and $L=K$)
this is the usual Milnor--Wood inequality (cf.\
\cite{bradlow-garcia-prada-gothen:2003}). 

The study of properties of Higgs bundles with extremal values for the
Toledo invariant is an interesting question. This has been studied for
various specific groups $G$ of hermitian type by Hitchin
\cite{hitchin:1987a} for $\PSL(2,\R)$, Gothen \cite{gothen:2001} for
$\Sp(4,\R)$, Garc{\'\i}a-Prada--Gothen--Mundet
\cite{garcia-gothen-mundet:2013} for $\Sp(2n,\R)$,
Bradlow--Garc{\'\i}a-Prada--Gothen
\cite{bradlow-garcia-prada-gothen:2003,bradlow-garcia-prada-gothen:2004,bradlow-garcia-gothen:2015}
for $\SO^*(2n)$ and $\U(p,q)$. A general study for $G$-Higgs bundles
for non-compact groups of hermitian type was carried out by
Biquard--Garc{\'\i}a-Prada--Rubio
\cite{biquard-garcia-rubio:2015}. From the point of view of
representations of surface groups much work has also been done and
without being at all exhaustive, we mention here a few works: Toledo
\cite{toledo:1979}, Hern\'andez \cite{hernandez:1991} and
Burger--Iozzi--Wienhard
\cite{burger-iozzi-wienhard:2003,burger-iozzi-wienhard:2010}. From
either point of view, one of the key properties of maximal objects
(Higgs bundles or representations) is that they exhibit rigidity
phenomena, of which we mention but two examples. Firstly, a
classical theorem of Toledo, which states that a maximal
representation of $\pi_1X$ in $\U(p,1)$ factors through
$\U(1,1)\times\U(p-1)$.  Secondly we mention
\cite[Proposition~3.30]{bradlow-garcia-prada-gothen:2003}, which says
that the moduli space of maximal $\U(p,p)$-Higgs bundles is isomorphic
to the moduli space of $K^2$-twisted Hitchin pairs of rank $p$ --- so
here Hitchin pairs play an important role even in the theory of usual
Higgs bundles. Toledo's theorem and its generalizations for surface
group representations have clear parallels on the Higgs bundle side of
the non-abelian Hodge theory correspondence. On the other hand, the
surface group representation parallel of the second kind of rigidity
phenomenon is perhaps less clear; see, however, Guichard--Wienhard
\cite{guichard-wienhard:2010} for the case of representations in
$\Sp(2n,\R)$.

\section{Wall crossing for $\U(p,q)$-Hitchin pairs}
\label{sec:wall-crossing}

We finish this paper by describing an application of wall-crossing
techniques to moduli of $\U(p,q)$-Hitchin pairs, following
\cite{nozad:2016,gothen-nozad:2016}. These techniques have a long
history in the subject, going back at least to Thaddeus' proof
\cite{thaddeus:1994} of the rank 2 Verlinde formula.  The main results
on connectedness of moduli spaces of $\U(p,q)$-Higgs bundles from
\cite{bradlow-garcia-prada-gothen:2003} were based on the
wall-crossing results for \emph{triples} of
\cite{bradlow-garcia-prada-gothen:2004}: triples are $Q$-bundles for a
quiver with two vertices and one arrow between them, so they
correspond to $\U(p,q)$-Hitchin pairs with one of the Higgs fields
$\beta$ or $\gamma$ vanishing. Later some of these results have been
generalized to \emph{holomorphic chains}, i.e., $Q$-bundles for a
quiver of type $A_n$, see
{\'A}lvarez-C{\'o}nsul--Garc{\'{\i}}a-Prada--Schmitt
\cite{alvarez-consul-garcia-prada-schmitt:2006},
Garc{\'{\i}}a-Prada--Heinloth--Schmitt
\cite{garcia-heinloth-schmitt:2014}, Garc{\'{\i}}a-Prada--Heinloth
\cite{garcia-heinloth:2013} and Heinloth \cite{heinloth:2014}. Similar
ideas have also been employed by other authors to study various
properties of moduli spaces, including their Hodge numbers, such as the
works of Bradlow--Garc{\'{\i}}a-Prada-- Mu{\~n}oz--Newstead
\cite{bradlow-garcia-munoz-newstead:2003},
Bradlow--Garc{\'{\i}}a-Prada--Mercat--Mu{\~n}oz--Newstead
\cite{bradlow-garcia-mercat-munoz-newstead:2007}, Mu\~noz
\cite{munoz:2008,munoz:2009,munoz:2010} and
Mu\~noz--Ortega--V{\'a}zquez-Gallo
\cite{munoz-ortega-vazques:2007,munoz-ortega-vazques:2009}.

One common feature of all these results is that they deal with quivers
without oriented cycles, corresponding to nilpotent Higgs fields. It
is therefore interesting to investigate to what extend the
aforementioned results can be generalized to quivers with oriented
cycles. Since we need at least two vertices to have effective stability
parameters, the simplest possible case is that of $\U(p,q)$-Hitchin
pairs, corresponding to the quiver \eqref{eq:upq-quiver}.

It turns out that a direct generalization of the arguments for triples
of \cite{bradlow-garcia-prada-gothen:2004} runs into difficulties. To
explain this, we first remark that the stability condition can only
change for certain discrete values of the parameter $\alpha$, called
\emph{critical values}. Fix topological invariants $t=(p,q,a,b)$ of
$\U(p,q)$-Hitchin pairs, where $a=\deg(V)$ and $b=\deg(W)$. Then
$\alpha$ is a critical value of the stability parameter for
$\U(p,q)$-Hitchin pairs of type $t$ if it is numerically possible to
have a proper subobject $E'\subset E$ of
a $\U(p,q)$-Hitchin pair $E=(V,W,\beta,\gamma)$ of type $t$ such that
\begin{equation}
  \label{eq:critical}
  \mu_{\alpha}(E') = \mu_{\alpha}(E)
  \quad\text{and}\quad
  \frac{p'}{p'+q'}\neq \frac{p}{p+q}
\end{equation}
(Here the type of $E'$ is $t'=(p',q',a',b')$.)
This means that $\alpha$ is critical if and only if it is possible for
$\U(p,q)$-Hitchin pairs to exist which are $\alpha'$-stable for
$\alpha'<\alpha$ and $\alpha'$-unstable for $\alpha'>\alpha$ (and
vice-versa). Denote by $\mathcal{M}_{\alpha^\pm}$ the moduli space
of $\alpha^\pm$-semistable $\U(p,q)$-Hitchin pairs of type $t$, where
$\alpha^\pm=\alpha\pm\epsilon$ for $\epsilon>0$ small. Then one is led
to introduce ``flip loci''
$\mathcal{S}_{\alpha^\pm}\subset\mathcal{M}_{\alpha^\pm}$
corresponding to $\U(p,q)$-Hitchin pairs which change their stability
properties as the critical value $\alpha$ is crossed. If one can
estimate appropriately the codimension of these flip loci, it will
follow that $\mathcal{M}_{\alpha^\pm}$ are birationally
equivalent. The $\U(p,q)$-Hitchin pairs $E$ in the flip loci have
descriptions as extensions
\begin{displaymath}
  0\to E'\to E \to E''\to 0
\end{displaymath}
for $\alpha$-semistable $\U(p,q)$-Hitchin pairs (of lower rank) $E'$
and $E''$ satisfying \eqref{eq:critical}. Such extensions are
controlled by the first hypercohomology of a two-term complex of
sheaves $\underline{Hom}^\bullet(E'',E')$ (see \cite[Definition
2.14]{gothen-nozad:2016}, cf.\ \cite{gothen-king:2005}). Thus, in
order to control the number of extensions one needs vanishing results
for the zeroth and second hypercohomology groups. This,
together with an analysis of the moduli space for large $\alpha$, was
the strategy followed in \cite{bradlow-garcia-prada-gothen:2003} to
prove irreducibility of moduli spaces of holomorphic triples. 

The main difficulty in generalizing this approach to $\U(p,q)$-Hitchin
pairs is that the vanishing results do not generalize without
additional hypotheses (compare, for example,
\cite[Proposition~3.6]{bradlow-garcia-prada-gothen:2004} and \cite
[Proposition 3.22]{gothen-nozad:2016}).  However, for a certain range
of the parameter $\alpha$ and the Toledo invariant, things can be made
to work.
Thus
we can obtain birationality of moduli spaces of $\U(p,q)$-Hitchin
pairs under certain constraints (see
\cite[Theorem~5.3]{gothen-nozad:2016}). This combined with the results
from \cite{bradlow-garcia-prada-gothen:2003} on connectedness of
moduli of $\U(p,q)$-Higgs bundles finally gives the main result: 

\begin{theorem}[{\cite[Theorem~5.5]{gothen-nozad:2016}}]
Denote by $\mathcal{M}_\alpha(p,q,a,b)$ the moduli
space of semistable $K$-twisted $\U(p,q)$-Hitchin pairs. Suppose that
$\tau=\frac{2pq}{p+q}(a/p-b/q)$ satisfies
$\abs{\tau}\leq\min\{p,q\}(2g-2)$. 
Suppose also that either one of the following conditions holds:
\begin{itemize}
\item[$(1)$] $a/p-b/q>-(2g-2) $, $q\leq p$ and $0\leq\alpha<\frac{2pq}{pq-q^2+p+q}\big(b/q-a/p-(2g-2)\big)+2g-2$,
\item[$(2)$] $a/p-b/q<2g-2$, $p\leq q$ and $\frac{2pq}{pq-p^2+p+q}(b/q-a/p+2g-2)-(2g-2)<\alpha\leq 0$.
\end{itemize}
Then the closure of the stable locus in the moduli space
$\mathcal{M}_\alpha(p,q,a,b)$ is irreducible.  In particular, if
$\gcd(p+q,a+b)=1$, then $\mathcal{M}_\alpha(p,q,a,b)$ is irreducible.
\end{theorem}

\begin{remark}
  Unless $p=q$, the conditions on $a/b-b/q$ in the theorem are
  guaranteed by the hypothesis $\abs{\tau}\leq\min\{p,q\}(2g-2)$
  (see~\cite[~Remark~5.6]{gothen-nozad:2016}). 
\end{remark}


\providecommand{\href}[2]{#2}

\noindent 
      Peter B.\ Gothen \\
      Centro de Matem\'atica da Universidade do Porto \\
      Faculdade de Ci\^encias da Universidade do Porto\\
      Rua do Campo Alegre, s/n \\
      4169--007 Porto -- Portugal\\
      \texttt{pbgothen@fc.up.pt}
\end{document}